\documentclass[11pt,reqno,final]{article}

\usepackage{type1cm}        

\usepackage{makeidx}         % allows index generation
\usepackage{graphicx}        % standard LaTeX graphics tool
\usepackage{multicol}        % used for the two-column index
\usepackage[bottom]{footmisc}% places footnotes at page bottom
\usepackage{a4wide, mathrsfs}
\usepackage{amsmath, amsfonts}

\renewcommand{\div}{\operatorname{div}}

 \usepackage{url}

\title{\bf On the segmentation of astronomical images via level-set methods}

\author{Silvia Tozza\footnotemark[1]  \and Maurizio Falcone\footnotemark[2]  }
\begin{document}
		
		\maketitle
		
	\renewcommand{\thefootnote}{\fnsymbol{footnote}}
	
	\footnotetext[1]
	{
	Istituto Nazionale di Alta Matematica, U.O. Dipartimento di Matematica,  ``Sapienza" Universit{\`a}  di Roma,
	P.le Aldo Moro, 5 - 00185 Rome, Italy
	({\tt e-mail: tozza@mat.uniroma1.it})
	}
	
	\footnotetext[2]
	{
		Dipartimento di Matematica,  ``Sapienza" Universit{\`a}  di Roma,
		P.le Aldo Moro, 5 - 00185 Rome, Italy
		({\tt e-mail: falcone@mat.uniroma1.it})
	}
	
	%\footnotetext[3]
	%{
\thanks{The authors are members of the INdAM Research group GNCS.}
	%}

\abstract{Astronomical images are of crucial importance for astronomers since they contain a lot of information about celestial bodies that can not be directly accessible. Most of the information available for the analysis of these objects starts with sky explorations via telescopes and satellites. Unfortunately, the quality of astronomical images is usually very low with respect to other real images and this is due to technical and physical features related to their acquisition process. This increases the percentage of noise and makes more difficult to use directly standard segmentation methods on the original image. 
In this work we will describe how to process astronomical images in two steps:  in the first step we improve the image quality  by a rescaling of light intensity whereas in the second step we apply level-set methods to identify the objects. Several experiments will show the effectiveness of this procedure and the results obtained via various discretization techniques for level-set equations.

%%%%%%%%%%%%%%%%%%%%%%%%%%%%%%%%%%%%%%%%%%%%%%%%%%%%%%%%%
\section{Introduction}\label{sec:intro}

Astronomical images  are acquired by appropriate sensors, called CCDs (\emph{Charge-Coupled Devices}), that are able to generate an electric charge at each pixel.  This charge is directly proportional to the electromagnetic radiation that affects the pixel and is the measure corresponding to the ``brightness'' of real optical images. A typical feature of astronomical images is that they suffer from various types of noise which make difficult to analyze them. Their noise percentage is usually much higher than that of standard optical images since the value at every pixel does not correspond to the flow of photons emitted from the light source, i.e.  the real signal, but  is modified by the  disturbances in the acquisition process. Let us recall the most important disturbances: 
\begin{itemize}
	\item the  noise related to the signal, modeled by a Poisson distribution with standard deviation $\sqrt{n_e}$, which is directly proportional to the flux emitted by the source
	\item the light coming from other celestial bodies and from the sky, i.e. the spurious light collected by the telescope (the so-called \emph{sky background})
	\item the thermal noise, caused by overheating of the CCD sensors, which leads to an increase of the thermal agitation and the generation of additional conduction electrons;
	\item the \emph{readout noise}, caused by the  electronic components of the CCD and due to the discrete nature of the signal.
\end{itemize}
The amount of noise present in the image is expressed mathematically in terms of  $SNR$ (\emph{Signal to Noise Ratio}), defined as the ratio between the power of the represented signal and that of the estimated noise, considering all the components previously listed. Larger values for this ratio correspond to images of better quality. 
The original image can not be used  for an accurate scientific analysis of the data as we will see in the following sections. For that reason, a series of preprocessing steps are performed to reduce the noise and improve the image quality. It has been shown that, by increasing the exposure time of the sensors to light, the ratio between signal and noise can be greatly increased. This improvement is directly proportional to $\sqrt{t_{exp}}$ but an exposure time that is too long can lead to a saturation of the pixels so this procedure has to be carefully implemented.  
Furthermore, noise reduction operations are performed on each image. Typical operations include  masking the defective pixels, subtracting the estimated value for the sky background and calibrating the image, but one can also apply a standard (linear or nonlinear) filter as we will do in our experiments.  After these operations the value of the flow, with its relative uncertainty, and the astronomical coordinates associated to each pixel  are redefined.  Among the many other precautions that can be used, we emphasize that  the most recent astronomical instruments use  cooling devices for the CCD sensors which allow to reduce the readout noise. Despite the operations of calibration and noise reduction and the wide variety of techniques that has been adopted, noise remains one of the main components of the astronomical images. 
Due to the above steps in the acquisition, the range for the admissible values for the astronomical images is really different from the range of other  kinds of images, e.g. it is common to have negative values at some pixels after the subtraction of the sky background. 
A final difference with respect to classical images is   the format currently used to store astronomical images. The most common format is FITS (\emph{Flexible Image Transport System}, \cite{NasaFITS}), introduced by the International Astronomical Union FITS Working Group (IAUFWG) in $1981$ and up-dated in $2016$ to its fourth version. The introduction of a new format is due to the need of save different information related to the images generated through  CCD sensors, such as the angular coordinates of the portion of sky observed or the zero-point magnitude of the sensor used.  This makes necessary to establish a common format, through which all the astronomers can interpret the data  in the same way. The format has been developed so that all files, even the oldest ones, can be read from every machine, structuring files as a sequence of logical data.\\
To set this paper into perspective, let us mention some related contributions in the literature. The problem of deblurring astronomic images produced by telescopes is a classical and difficult problem in the astronomical community \cite{F65, WG91, R04}. A novel technique to reduce the distortion caused by the ground-level turbolence of the atmosphere has been recently proposed in \cite{MLSB14}. A similar goal has motivated the development of a high-resolution speckle imaging technique presented  in \cite{HJHN16}. As far as segmentation models is concerned, we mention that a modified version of the Chan-Vese model \cite{CV01} has been proposed and analyzed in \cite{GB11}, some results obtained by a high-order splitting scheme are also presented there. It is interesting to note that this is  a region based method with a level-set representation that can be applied to multispectral images. 

In this paper we propose a strategy to analyze and segment astronomical images via the level-set method introduced in \cite{OS88}. Although the segmentation problem has been investigated by many authors (see e.g the monographies \cite{S99, OF03, FF13} and the references therein) and several successful applications have been reported in many areas,  level-set techniques are still not very popular in the astronomers community. Most probably this  is due to the above mentioned features of astronomical images that make a direct application of these methods fail or give inaccurate results. Here we propose a coupling between an appropriate rescaling technique and a standard level-set methods to improve the global accuracy of the segmentation and increase the number of celestial bodies that can be extracted from a single image. We also add a  filtering step  to reduce the noise before segmenting. Hopefully, this will help astronomers in their sky investigations. \\
The paper is organized as follows: In Sect. \ref{sec:astro-images}, we propose new different rescaling transformations, adopted as the first two steps of our algorithm to improve the results of the segmentation of astronomical images.  We briefly recall in Sect. \ref{sec:segment_LS} the first and second order equations related to level-set methods and the corresponding finite difference and semi-Lagrangian schemes that we used for our numerical experiments, at the beginning of this section we give some hints on the filtering step. Finally, in Sect. \ref{sec:tests}, we present  our complete Rescaling Segmentation Algorithm (RSA) and we discuss in detail our numerical tests on simulated and real astronomical images. We conclude with Sect. \ref{sec:conclusions} where we summarize our  final remarks and future perspectives. 

%%%%%%%%%%%%%%%%%%%%%%%%%%%%%%%%%%%%%%%%%%%%%%%%%%%%%%%%%%%%%
\section{Efficient rescaling of astronomical images} \label{sec:astro-images}
Let us start describing the first step of the procedure we adopted to segment astronomical images.  It is useful to read astronomical images saved in the FITS format in MATLAB, thanks to the command \emph{fitsread} and transform them in the matrix format that is common in image processing. The matrix $I_0$ returned as output from the function \emph{fitsread} can take negative values due to the preprocessing techniques of calibration and reduction of noise applied to the images provided by the CCD sensors (e.g. procedures as the calibration or the subtraction operation of the estimated sky background).

\noindent We need to rescale the image values, defined on a rectangular domain $\overline{\Omega}$, with $\Omega\subset \mathbb{R}^2$,  in order to obtain real values in $[0,1]$. Starting from the matrix $I_0$, this is done defining 
\begin{equation}\label{risc01}
\widetilde{I_0}=\frac{I_0(x,y)-m_0}{M_0-m_0}\,,
\end{equation}
where
\[
m_0:=\min_{(x,y)\in\overline{\Omega}}I_0(x,y)\,,\quad M_0:=\max_{(x,y)\in\overline{\Omega}}I_0(x,y)\, .
\]
The image $\widetilde{I_0}$ obtained is still not  ready for the segmentation since, in most cases, is very dark and only few celestial bodies will be visible to the naked eye. For that reason, we choose to transform the image, rescaling the values of the pixels by means of an appropriate  function that we will construct in the sequel.

%%%%%%%%%%%%%%%%%%%%%%%%%%%%%%%%%%%%%%%%%%%%%%%%%%%%%%%%%%%%%%%%%%%%%%%%%%%%%%%%%%%
\subsection{Elevation to power or logarithmic rescaling}\label{sectrisc1}
We look for a rescaling  function $r:[0,1]\rightarrow\mathbb{R}$ for the gray levels. Since these values for the image $\widetilde{I_0}$ obtained by  \eqref{risc01} are in the range $[0,1]$, the function $r$ must satisfy the following conditions:
\begin{enumerate}
	\item [ A1. ] $r([0,1])\subseteq[0,1]$ \label{pr1}
	\item [ A2. ]$r(0)=0$, $r(1)=1$ \label{pr2}
	\item[ A3. ] $r$ strictly increasing. \label{pr3}
\end{enumerate}
In other words, the rescaling transformation must keep the minimum and maximum brightness points of the image unaltered and rescale the intermediate values, without changing their ordering. Since the image is very dark, we also want the transformation to amplify the brightness values. In mathematical terms, we require $r$ to satisfy the additional condition
\begin{enumerate}
	\item[A4. ] $r(x)>x\,,\quad\forall x\in[0,1]$.
\end{enumerate}
Clearly, several functions can satisfy the above four properties. A simple choice is given by
\begin{equation}
r_1(x):=x^\alpha\,,
\end{equation}
with $\alpha\in(0,1)$ a fixed parameter. \\
Another function can be obtained by a logarithmic transformation of the form
\begin{equation}
r_2(x):=\left[\frac{\ln(x+1)}{\ln 2}\right]^\alpha\,,
\end{equation}
with $\alpha\in(0,1)$. 
Both functions converge pointwise to the identity for $\alpha$ going to 1, whereas for $\alpha$ going to 0 they converge pointwise to the function
	\begin{equation}
	\widetilde{r_1}(x):=
	\begin{cases}
	0\,\quad&\text{if }x=0\,,\\
	1\, &\text{if }x\in(0,1]\,.
	\end{cases}
	\end{equation}
The latter transforms every brightness value, with the exception of the null one, assigning it the value $1$. 
After the rescaling, the brightness increases as $\alpha$ decreases. %The behavior of the two functions, as the $\alpha$ changes, is visible in Fig. \ref{funzrisc1}. 
The behavior of the two functions for different values of $\alpha$  is visible in Fig. \ref{funzrisc1}. \\
\begin{figure}[h]%[!tb]
	\centering
	\includegraphics[width=0.48\textwidth]{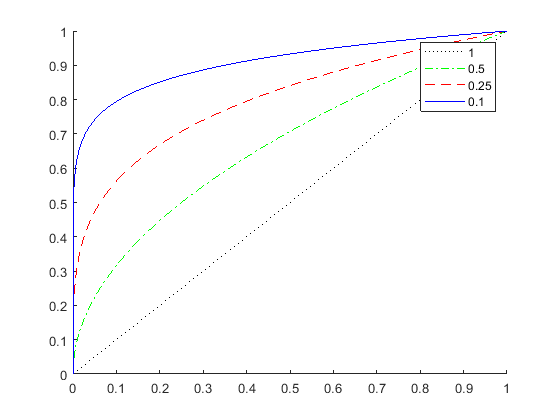}
	%\quad
	\includegraphics[width=0.48
	\textwidth]{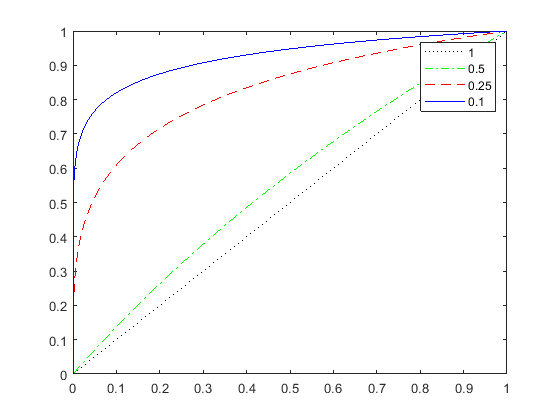}
	\caption{Performance of the functions $r_1$ and $r_2$, by varying the parameter $\alpha\in(0,1]$.}
	\label{funzrisc1}
\end{figure}
In the numerical tests presented in Section \ref{sec:tests}, %\ref{TestAstro}, 
we will only show the results obtained with the function $r_1$, since for the same $\alpha$,  $ r_2 $ gives  almost identical results.

%%%%%%%%%%%%%%%%%%%%%%%%%%%%%%%%%%%%%%%%%%%%%%%%%%%%%%%%%%%%%%%%%%%%%%%%%%%%%%%%%%%%%
\subsection{Rescaling with a threshold}\label{sectrisc2}
From our experiments on astronomical images (see Sect. \ref{sec:tests})  we have observed that the proposed transformations $r_1$ and $r_2$  can be improved. As we said,  astronomical images are affected by a strong  noise component and the rescaling  has a significant effect also on  high brightness values due to the noise component. When these values are rescaled, they result too high so the global effect is an amplification of the tone differences with respect to the pixels closer to the real tone of the background. This amplification can make the segmentation method fail, identifying artificial objects that are not present in the real image. 
To avoid this undesired effect  the rescaling should distinguish the pixels of celestial bodies from those of the background: the values of the former must be amplified, while the others must be attenuated. A natural idea is to introduce a threshold to determine the gray tones of the objects and, to be optimal,  this threshold should be automatically identified by an algorithm. A good choice for standard images  is provided by the Otsu algorithm \cite{Otsu79}, so we decided to use the value $\tau \in [0,1) $ provided by this algorithm as a threshold. 
The new rescaling transformation must still respect the properties A1--A3 of Subsect. \ref{sectrisc1}, in addition it has to satisfy the condition 
\begin{equation}\label{pr4soglia}\tag{A4$_{\tau}$}
	\begin{cases}
	r(x)<x\,,\quad&\text{if }0<x<\tau\,,\\
	r(\tau)=\tau\,,\\
	r(x)>x\,, &\text{if }\tau<x<1\,,
	\end{cases}
\end{equation}
and to be continuous at $x=\tau$. \\
A rescaling function of this type can be obtained by considering the applications $x^\beta$ and $x^{1/\beta}$, with $\beta\in\mathbb{N}\setminus\{0\}$, 
respectively in the two subsets $[0, \tau]$ and $(\tau, 1]$. These functions have to be appropriately translated and expanded to respect all the conditions.  In this way, we obtain the function 
\begin{equation}
	r_3(x):=\begin{cases}
	\displaystyle\frac{x^{\beta}}{\tau^{\beta-1}}\,,\quad&\text{if }0\leq x<\tau\,,\\[1.5ex]
	\displaystyle\frac{(x-\tau)^{1/\beta}}{(1-\tau)^{1/\beta-1}}+\tau\,,\quad &\text{if }\tau\leq x\leq 1\,.
	\end{cases}
\end{equation}
The behavior of $r_3$ varying $\beta\in\mathbb{N}\setminus\{0\}$ is reported in Fig. \ref{funzrisc2}. 
This function satisfies the properties listed before, converges pointwise to the identity function for $\beta$ tending to 1 and to the following function 
\begin{equation}
	\widetilde{r_3}(x):=\begin{cases}
	0\,,\quad&\text{if }0\leq x<\tau\,,\\
	\tau\,, &\text{if }x=\tau\,,\\
	1\,, &\text{if }\tau<x\leq 1\,,
	\end{cases}
\end{equation}
for  $\beta$ tending to  $+\infty$. 
\begin{figure}[h]%[!tb]
		\centering
		\includegraphics[width=0.68\textwidth]{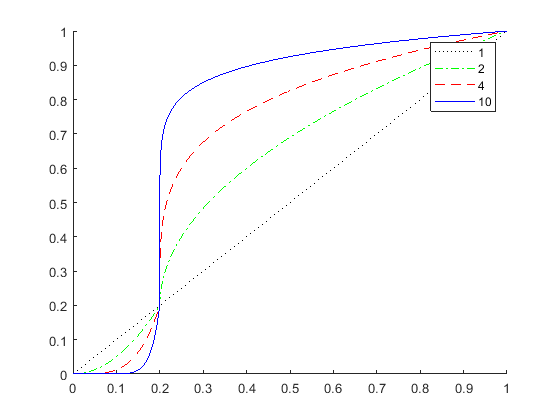}
		\caption{Behavior of the function $r_3$ by varying the parameter $\beta\in\mathbb{N}\setminus\{0\}$.}
		\label{funzrisc2}
\end{figure}
%%%%%%%%%%%%%%%%%%%%%%%%%%%%%%%%%%%%%%%%%%%%%%%%%%%%%%%%%%%%%%%%%%
%%%%%%%%%%%%%                    SEGMENTATION VIA LEVEL SET METHODS
%%%%%%%%%%%%%%%%%%%%%%%%%%%%%%%%%%%%%%%%%%%%%%%%%%%%%%%%%%%%%%%%%%
\section{Segmentation  via level-set methods}\label{sec:segment_LS}
As we said in the introduction, we  follow the level-set (LS) approach to segmentation problems obtained by the rescaling procedure described in the previous section. For readers convenience let us briefly describe the main features of this approach. The level-set method has been introduced by Osher and Sethian \cite{S85,OS88} and since then it has been  widely used in many  applications, e.g.  fronts propagation, computer vision, computational fluids dynamics (see \cite{S99,OF03} for several interesting examples). Its popularity is due to the simplicity in the implementation and its capability to follow topological changes (splitting, merging)  in time. Typical examples are when a planar curve (or a multidimensional surface) splits into many parts or when several evolving curves (or surfaces)  merge into a single one. This is the main reason for its popularity also in the image processing community. For the segmentation problem the idea is to define a normal vector field bringing an initial curve (e.g. a circle) onto the object boundaries in an image. 

Let us consider an image $\widetilde{I_0}:\overline{\Omega}\rightarrow[0,1]$, with $\Omega\subset\mathbb{R}^2$ an open rectangular domain. 
Let us fix an initial curve $\gamma_0\subset\overline{\Omega}$. We want to track its evolution according to the normal velocity  and we define it  so that it goes to zero (and therefore the front stops) at the edges of the object to be identified. The methods based on this approach can be divided into two subclasses. In the methods belonging to the  first class,  the speed depends on the gradient of the image $\widetilde{I_0}$ at each point $(x, y) \in \Omega$, since the gradient provides a measure of the gray-level variation in the image and therefore it identifies the presence of edges. The second class of methods, introduced by Chan and Vese \cite{CV01}, is inspired by a variational segmentation technique proposed by Mumford and Shah \cite{MS89} and is based on the minimization of a functional which allows to partition the image in  regions  where there is a small variation of gray levels.  
Looking more in details the first class, we have to solve an evolutive Hamilton-Jacobi equation
\begin{equation}
\label{levelsetd2}
\begin{cases}
u_t(x,y,t)+c(x,y,t)\lvert\nabla u(x,y,t)\rvert=0\,,\quad &\forall(x,y,t)\in\Omega\times(0,T]\,,\\
u(x,y,0)=u_0(x,y)\,,&\forall(x,y)\in\overline{\Omega}\,,
\end{cases}
\end{equation}
with $u(\cdot,\cdot,t)$ and  $u_0$ the representation function of the front at time $t$ and at the initial time, respectively, and $c$ is the velocity function. 
Depending on the definition for $c$ (that in general may depend on $x$, $t$ and the curvature), Eq. \eqref{levelsetd2} will be a first or second order equation. Several explicit definitions of $c$ will be reported in Sect. \ref{subsec:edge-detector}.
In order to segment a given image $\widetilde{I_0}$, we choose the initial front $\gamma_0 \subset \overline{\Omega}$ and its representation function $u_0: \overline{\Omega} \rightarrow \mathbb{R}$. In particular, if we want to approximate the edges of the object with a curve that expands from within, we choose $u_0$ in such a way that, denoted by $\omega_0$ the region of the plane enclosed by the front $\gamma_0$, with $\gamma_0=\partial\omega_0$ and $\omega_0$ open, we put
\begin{equation}
\begin{cases}
u_0(x,y)<0\,,\quad&\forall(x,y)\in\omega_0\,,\\
u_0(x,y)=0\,, &\forall(x,y)\in\gamma_0\,,\\
u_0(x,y)>0\,, &\forall(x,y)\in\overline{\Omega}\setminus\overline{\omega_0}\,.
\end{cases}
\end{equation}
Conversely,  if we want the front to contract, we can reverse  the sign of the  initial representative or of the normal direction. 
Next, we fix the velocity of the front and we solve the equation of the level-set method obtained by it: denoted by $u$ its solution, we obtain the front at time $t> 0$ as the $0$-level-set of $u(\cdot,\cdot,t)$, that is
\begin{equation}
\gamma_t=\left\{(x,y)\in\overline{\Omega}\,|\,u(x,y,t)=0\right\}\,.
\end{equation}
Equation \eqref{levelsetd2} is complemented with boundary conditions. We chose to use homogeneous Neumann conditions
\begin{equation}
\frac{\partial u}{\partial\mathbf{\eta}}(x,y,t)=0\,, \quad\forall(x,y,t)\in\partial\Omega\times(0,T]\,.
\end{equation}
The choice of the final time $T$ to which numerically solve the equation \eqref{levelsetd2} will have to be carried out through a stopping criterion, which detects when the front is near equilibrium, through the verification of a condition. In this paper, we will adopt the following criterion:
First, at each iteration we identify the grid nodes near the front with respect to a fixed tolerance denoted by $\varepsilon_F$. More precisely, since the front at time $t_n$ is the $0$-level curve of the representation function, we define the approximate front by means of $\mathbf{V}^n := \{v_ {i, j}^n \}$, where  $v_{i,j}^n$ is the value computed on the grid node $(x_{i},y_j)$ at time $n$. Our numerical front is given by 
\begin{equation}
F^n \equiv \{ (x_i, y_j) : \lvert v_{i,j}^n\rvert\leq\varepsilon_F\,\}.
\end{equation}
Let us denote by $\mathscr{F}$ the set of indexes of the nodes that respect this condition and with $\mathbf{V}^{n,F}$ the vector formed by the elements of $\mathbf{V}^n$  corresponding to them. Hence, we fix an additional tolerance $\varepsilon$: the stopping condition of the numerical scheme will be
\begin{equation}\label{critarr}
\left\|\mathbf{V}^{n+1,F}-\mathbf{V}^{n,F}\right\|_1\leq \varepsilon\,,
\end{equation}
with the norm $\|\cdot\|_1$ defined by 
\begin{equation}
\left\|\mathbf{V}^{n+1,F}-\mathbf{V}^{n,F}\right\|_1:=\Delta x^2\sum_{(i,j)\in\mathscr{F}} \left\lvert v_{i,j}^{n+1}-v_{i,j}^n\right\rvert\,.
\end{equation}
In other words, we proceed to solve the scheme up to the  $(n+1)$-th iteration when the representation has reached equilibrium with a tolerance  $\varepsilon$  at all the nodes belonging to $F^n$.

%%%%%%%%%%%%%%%%%%%%%%%%%%%%%%%%%%%%%%%%%%%%%%%%%%%%%%%%%%%%%%%%%%%%%%%%%%%%%%
\subsection{The filtering pre-processing step}\label{subsec:filter}
Let us analyze the first class of active contours methods. Since the edges of objects are, in most cases, identified by large variations of gray tones in their neighborhood, we can define the velocity of the front as a function of the gradient of the function $\widetilde{I_0}$ that models the image. However, $\widetilde{I_0}$ is a noisy image so in order to define its gradient it is useful to add a filtering step on it. We did it in two different ways: by applying a Gaussian filter, i.e. solving the \emph{heat equation with homogeneous Neumann conditions}
\begin{equation}\label{EqCalore}
\begin{cases}
I_t(x,y,t)=\Delta I(x,y,t)\,,\quad&\forall(x,y,t)\in\Omega\times(0,T_C]\,,\\[1.5ex]
\displaystyle\frac{\partial I}{\partial\mathbf{\eta}}(x,y,t)=0\,, &\forall(x,y,t)\in\partial\Omega\times(0,T_C]\,,\\[1.5ex]
I(x,y,0)=\widetilde{I_0}(x,y)\,,&\forall(x,y)\in\overline{\Omega}\,,
\end{cases}
\end{equation}
which has a diffusive effect on the initial datum $\widetilde{I_0}$, for a small fixed time $T_c > 0$ (in the numerical tests, it is of the order of $10^{-3}$ or $10^{-4}$). 
Numerically, we solve \eqref{EqCalore} by the standard centered finite difference  scheme
\begin{equation}\label{schemacal}
I_{i,j}^{n+1}=I_{i,j}^n+\widetilde{\Delta t}\left[\frac{I_{i+1,j}^n+I_{i,j+1}^n-4I_{i,j}^n+I_{i-1,j}^n+I_{i,j-1}^n}{\Delta x^2}\right]\,,
\end{equation}
forward in time, with time step  $\widetilde{\Delta t} >0$ and space steps $\Delta x= \Delta y$.  In \eqref{schemacal} $I_{i,j}^n$ denotes as usual the approximation of the gray level of the image at  the pixel of coordinate $(i, j)$ and at time $t_n$, whereas $I^0_{i,j}:=\widetilde{I_0}(x_i,y_j)$ for every $(i,j)\in\mathscr{I}$, set of indexes. 
%sottoinsieme di Z2 tale che $(x_i,y_j) \in \Omega$ per ogni $(i,j)\in\mathscr{I}$.
The required CFL condition for this numerical scheme is $\widetilde{\Delta t}\leq\Delta x^2/4$.  \\
The consequence of applying the Gaussian filter is an
edge blurring due to isotropic diffusion. Choosing large values of $|\nabla I|$ as an indicator 
of the edge points of the image, we would like  to stop the diffusion at the edges, we pass  from an isotropic to an anisotropic diffusion, i.e.
\begin{equation}\label{eq:PM}
I_t %\frac{\partial I}{\partial t} 
= div(\nabla I) \hbox{ is replaced by }
%\frac{\partial I}{\partial t} 
I_t = div(f(|\nabla I|) \nabla I). 
\end{equation}
This is the idea behind the Perona-Malik model \cite{PM90} described by  \eqref{eq:PM} and complemented by suitable boundary conditions (e.g. homogeneous Neumann boundary conditions), the initial condition is the original image. The anisotropic diffusion is driven by $f$ and two typical choices for the diffusion coefficient are: % (also called edge-stopping functions):
\begin{eqnarray}
f_1(|\nabla I|) &= \displaystyle \frac{1}{1 + \Big(\frac{|\nabla I|}{\mu}\Big)^2}\\
f_2(|\nabla I|) &= \exp\Big(-\Big(\frac{|\nabla I|}{\mu}\Big)^2\Big)
\end{eqnarray}
where $\mu$ is the gradient magnitude threshold parameter. In our numerical simulations, we will use the function $f_2$. 
Let us denote by $\widetilde{I}_{filt}$ the solution of the problem \eqref{EqCalore} or \eqref{eq:PM}, filtered version of the  image $\widetilde{I_0}$. 

\subsection{Edge-detector functions}\label{subsec:edge-detector}
We want the  velocity $c$ of the front to vanish near the edges so we introduce a function $g$ of  $|\nabla \widetilde{I}_{filt}|$, called \emph{edge-detector}, that has to safisfy the following conditions:    
\begin{equation}
g:[0,+\infty)\rightarrow[0,+\infty) \hbox{ is decreasing and } \lim_{z\rightarrow+\infty}g(z)=0\,.
\end{equation}
In this way $g(\lvert\nabla \widetilde{I}_{filt}(x,y)\rvert)$ will tend to 0  approaching the points $(x, y)$ near the edges to be identified, since at the edges we typically have very  high values of $\lvert\nabla \widetilde{I}_{filt}\rvert$. Higher values of $g$ will correspond to points where $\lvert\nabla \widetilde{I}_{filt}\rvert\approx 0$, i.e. to the regions where the gray tones of the image are approximately constant. 
Two possible choices for the edge-detector function are the following: 
\begin{equation}\label{g1}
g_1(z):=\frac{1}{1+z^p}\,,\quad\forall z\in[0,+\infty)\,,  p\geq 1, 
\end{equation}
proposed in \cite{CCCD93} with $p = 2$, and in  \cite{MSV95} with $p=1$, 
and 
\begin{equation}\label{g2}
g_2(z):=1-\frac{z-m}{M-m}\,,\quad\forall z\in[0,+\infty)\,,
\end{equation}
where
\[
m:=\min_{(x,y)\in\Omega}\lvert\nabla \widetilde{I}_{filt}(x,y)\rvert \,,\quad M:=\max_{(x,y)\in\Omega}\lvert\nabla \widetilde{I}_{filt}(x,y)\rvert %\,.
\]
defined in \cite{MSV95}. 
Practically, the values of $g_1(\lvert\nabla \widetilde{I}_{filt}(x,y)\rvert)$ vary between $(1+M)^{-1}$ and $(1+m)^{-1}$, whereas the values of $g_2(\lvert\nabla \widetilde{I}_{filt}(x,y)\rvert)$ are between $0$ (for $\lvert\nabla \widetilde{I}_{filt}\rvert=M$) and $1$ (for $\lvert\nabla \widetilde{I}_{filt}\rvert=m$).\\
Let us discuss some typical choices for   the velocity $c$. A simple choice is to make $c$ dependent just on the point
\begin{equation}
c(x,y,t):=g(x,y)\,,\quad\forall(x,y)\in\Omega\,.
\end{equation}
In this way, using the notation $g(x,y):=g(\lvert\nabla \widetilde{I}_{filt}(x,y)\rvert)$, the problem to solve becomes
\begin{equation}\label{schtype1}
\begin{cases}
u_t(x,y,t)+g(x,y)\lvert\nabla u(x,y,t)\rvert=0\,,\quad &\forall(x,y)\in\Omega\,,\,\forall t\in(0,T]\,,\\
\frac{\partial u}{\partial\mathbf{\eta}}(x,y,t)=0\,, &\forall(x,y)\in\partial\Omega\,,\,\forall t\in(0,T]\,,\\
u(x,y,0)=u_0(x,y)\,,&\forall(x,y)\in\overline{\Omega}\,,
\end{cases}
\end{equation}
with $u_0$ the representation function of the initial front. 
This is the isotropic case and the corresponding equation is a first-order Hamilton-Jacobi equation of eikonal type.\\
Another popular choice is to use a  velocity that, at each point $(x, y)$, depends on the geometric properties of the front, e.g.  its curvature $k(x, y)$.  This choice is more complicated since the velocity will also depend on $u$, hence on $t$.  Following  \cite{ALM92,MSV93}, we consider a {\em curvature dependent velocity}
\begin{equation}
c(x,y,t):=g(x,y)\left(1-\nu k(x,y)\right)\,, \quad\forall(x,y)\in\Omega\,,
\end{equation}
where $\nu>0$ is a fixed parameter. 
The factor $g(x,y)$ causes that the front  stops near the edges. The parameter $\nu$ (typically less than $1$) weighs the speed dependency on the curvature. 
Since the curvature is given by   
\begin{equation}
k(x,y)=\div\left(\frac{\nabla u(x,u,t)}{\left\lvert\nabla u(x,y,t)\right\rvert}\right)\,,\quad\forall(x,y)\in\Omega\,,
\end{equation}
the level-set corresponding equation is the second order Hamilton-Jacobi equation
\begin{equation}
\label{schtype2}
\begin{cases}
u_t(x,y,t)+g(x,y)\lvert\nabla u(x,y,t)\rvert 
=\nu &\hspace{-0.3cm}g(x,y)\div\left(\frac{\nabla u(x,u,t)}{\left\lvert\nabla u(x,y,t)\right\rvert}\right) \lvert\nabla u(x,y,t)\rvert \,,\\
&\forall(x,y)\in\Omega\,,\,\forall t\in(0,T]\,,\\
\frac{\partial u}{\partial\mathbf{\eta}}(x,y,t)=0\,, &\forall(x,y)\in\partial\Omega\,,\,\forall t\in(0,T]\,,\\
u(x,y,0)=u_0(x,y)\,,  &\forall(x,y)\in\overline{\Omega}\,,
\end{cases}
\end{equation}
with the same buondary conditions and initial datum as in \eqref{schtype1}. 
The term in the second member has a diffusive effect on the solution: consequently, this type of scheme can be useful for segmenting images characterized by noise. 
Note that, in practice, the function $g$ is not necessarily equal to $0$ at all points on the edges of the objects, even if it takes very small values. The stopping  rule  \eqref{critarr} allows to control the numerical scheme so that  the evolution stops at time $T$ whenever the velocity is below a given threshold. \\
%%Schemi
In order to get a numerical solution of  \eqref{schtype1} and \eqref{schtype2}  in our tests we will use  a finite difference scheme (FD) and a semi-Lagrangian scheme (SL), so we will be able to compare their results.  Let us recall that the  {\em FD schemes}  for the two equations are, respectively, 
\begin{equation}\label{DFg}
v_{i,j}^{n+1}=v_{i,j}^n-\Delta t g_{i,j}\nabla^+\,,
\end{equation}
and
\begin{equation}
\label{DFgk}
\begin{cases}
v_{i,j}^{n+1}=v_{i,j}^n-\Delta t g_{i,j}\nabla^+ +\frac{\nu}{4} g_{i,j}(v_{i+1,j}^n 
 +v_{i,j+1}^n +v_{i-1,j}^n+v_{i,j-1}^n)\,\\
\hspace{5.8cm} \text{if }\lvert D_{i,j}^c[\mathbf{V}^n]\rvert\leq C\Delta x^s\,,\\[1.5ex]
v_{i,j}^{n+1}=v_{i,j}^n-\Delta t g_{i,j}\nabla^+ +\nu\Delta t g_{i,j}\Lambda_{(i,j)}^n\,,\quad \text{if }\lvert D_{i,j}^c[\mathbf{V}^n]\rvert> C\Delta x^s\,,
\end{cases}
\end{equation}
for each $(i,j)\in\mathscr{I}$ and $n\in\{0,1,\dots,N_T-1\}$, where $\mathbf{V}^0:=\mathbf{U}(0)$, $g_{i,j}:=g(\mathbf{x}_{i,j})$, with $\mathbf{x}_{i,j} := (x_i,y_j)$,
\begin{equation}
\begin{split}
\nabla^+&:=\big[\max\{D^{-x}_{i,j}[\mathbf{V}^n],0\}^2+\min\{D^{+x}_{i,j}[\mathbf{V}^n],0\}^2\\
&\quad+\max\{D^{-y}_{i,j}[\mathbf{V}^n],0\}^2+\min\{D^{+y}_{i,j}[\mathbf{V}^n],0\}^2\big]^{1/2}
\end{split}
\end{equation}
and
\begin{equation}
\begin{split}
\Lambda_{i,j}^n&:=\frac{1}{(D^{c,x}_{i,j}[\mathbf{V}^n]^2 +D^{c,y}_{i,j}[\mathbf{V}^n]^2)^{1/2}}\,\Big(D^{2,x}_{i,j}[\mathbf{V}^n] D^{c,y}_{i,j}[\mathbf{V}^n]^2\\
&\quad-2D^{c,x}_{i,j}[\mathbf{V}^n]D^{c,y}_{i,j}[\mathbf{V}^n]D^{xy}_{i,j}[\mathbf{V}^n] +D^{2,y}_{i,j}[\mathbf{V}^n]D^{c,x}_{i,j}[\mathbf{V}^n]^2\Big)\,.
\end{split}
\end{equation}
We refer the reader interested in the construction and the analysis of these schemes to \cite{S99, OF03}.\\
Let us also recall that the  {\em SL schemes} are, respectively,
\begin{equation}\label{SLg}
v_{i,j}^{n+1}=\min_{\mathbf{a}\in B(0,1)}\left\{I[\mathbf{V}^n](\mathbf{x}_{i,j}-\Delta t g_{i,j}\mathbf{a})\right\}
\end{equation}
and
\begin{equation}
\label{SLgk}
\begin{cases}
v_{i,j}^{n+1}=\min\limits_{\mathbf{a}\in B(0,1)}\left\{I[\mathbf{V}^n](\mathbf{x}_{i,j}-\Delta t g_{i,j}\mathbf{a})\right\} 
+\frac{\nu}{4} g_{i,j}(v_{i+1,j}^n+v_{i,j+1}^n+ v_{i-1,j}^n+v_{i,j-1}^n)\,, \\
\hspace{6cm} \text{if }\lvert D_{i,j}^c[\mathbf{V}^n]\rvert\leq C\Delta x^s\,,\\
v_{i,j}^{n+1}=\min\limits_{\mathbf{a}\in B(0,1)}\left\{I[\mathbf{V}^n](\mathbf{x}_{i,j}-\Delta t g_{i,j}\mathbf{a})\right\}
+\frac{\nu}{2}g_{i,j}\big[I[\mathbf{V}^n](\mathbf{x}_{i,j}+\sigma_{i,j}^n\sqrt{\Delta t})\\
\qquad\qquad+I[\mathbf{V}^n](\mathbf{x}_{i,j}-\sigma_{i,j}^n\sqrt{\Delta t})\big]\,,
 \hspace{1cm}\text{if }\lvert D_{i,j}^c[\mathbf{V}^n]\rvert> C\Delta x^s\,,
 \end{cases}
\end{equation}
for each $(i,j)\in\mathscr{I}$ and $n\in\{0,1,\dots,N_T-1\}$, with $\mathbf{V}^0:=\mathbf{U}(0)$ and
\begin{equation}
\sigma_{i,j}^n:=\frac{\sqrt{2}}{\lvert D_{i,j}^c[\mathbf{V}^n]\rvert}\left(\begin{matrix}D_{i,j}^{c,y}[\mathbf{V}^n]\\ -D_{i,j}^{c,x}[\mathbf{V}^n]\end{matrix}\right)\,.
\end{equation}
The role of the threshold $\Delta x^s$ for the first derivatives in \eqref{SLgk} is explained in detail in \cite{CFF10} (see also \cite{FF13} for the general theory of semi-Lagrangian schemes and \cite{CFF13} for other  applications to image processing problems). For our purposes it is sufficient to note that this threshold is used to solve also the degenerate case without adding a regularization.

%%%%%%%%%%%%%%%%%%%%%%%%%%%%%%%%%%%%%%%%%%%%
%%%%        SECTION NUMERICAL TESTS
%%%%%%%%%%%%%%%%%%%%%%%%%%%%%%%%%%%%%%%%%%%%
\section{Numerical Tests}\label{sec:tests}
Let us start describing the complete segmentation algorithm that includes the rescaling preprocessing via the functions $r_1$, $r_2$ and $r_3$ illustrated in Sec. \ref{sectrisc1} and \ref{sectrisc2}. \\

\noindent {\bf RESCALING SEGMENTATION ALGORITHM (RSA)} \\
\noindent {\em STEP 1:} Apply to the original image $I_0$ the rescaling defined in \eqref{risc01} to get $\widetilde{I_0}$ which takes values in $[0,1]$. \\
{\em STEP 2:} Choose one of the proposed rescaling functions $r_i$, $i\in\{1,2,3\}$, and set
\begin{equation}\label{trasff}
\widetilde{I}:=r_i(\widetilde{I_0})\,, 
\end{equation}
for each element of the matrix. 
In case we choose the function $r_3$, we first apply the thresholding method of Otsu to the image $\widetilde{I_0}$ in order to select the optimal threshold $\tau$, and then we apply  \eqref{trasff}. \\
{\em STEP 3: } Filter the image $\widetilde{I}$ by few iterations of the linear filter given by the scheme \eqref{schemacal}, or applying the PM  method \eqref{eq:PM} with $f_2$. This step produces $\widetilde{I}_{filt}$.  \\
{\em STEP 4: } Apply one of the segmentation active contours methods to the image $\widetilde{I}_{filt}$ thus obtained in  STEP 3.\\

%%%%%%%%%%%%%%%%%%%%%%%%%%%%%%%%%%%%%%%%%%%%%%%%%%%%%%%%%%%%%%%%%%
\vspace{0.2cm}
\noindent We are now ready to present some numerical tests, using the RSA algorithm. 
Let us consider an $M\times N$ image and let us fix the discretization parameters as: 
\begin{itemize}
\item $\Delta x=\Delta y=0.1$ the space step of the uniform grid 
\item $\Delta t =\Delta x/4=0.025$, for the  FD scheme approximating the first order problem
\item $\Delta t =\Delta x^2=0.01$, for the  FD scheme approximating the second order problem
\item $\Delta t =\Delta x=0.1$, for the SL schemes. 
\end{itemize}
That structured grid has nodes located at the center of the pixels, with coordinates $\big((j-1)\Delta x,(i-1)\Delta x\big)$, for $j=1,\dots,M$ and $i=1,\dots,N$, and the rectangular domain is defined as
\begin{equation}
\Omega:=\left[-\frac{\Delta x}{2},a-\frac{\Delta x}{2}\right]\times \left[-\frac{\Delta x}{2},b-\frac{\Delta x}{2}\right]\,,
\end{equation}
with  $a:=M\Delta x$ and $b:=N\Delta x$. 
For each test, we will consider three cases:
\begin{itemize}
\item a segmentation of the original image, without rescaling  (i.e. dropping STEP 2 of the algorithm, setting $\widetilde{I} = \widetilde{I}_0$)
\item a segmentation using a rescaling of the image by  $r_1$
\item a segmentation using a rescaling of the image by $r_3$ and the optimal threshold computed by the Otsu's algorithm.
\end{itemize}

\noindent As we said in Sect. \ref{sectrisc1}, we will omit the results obtained by  rescaling via $r_2$ since the results are almost identical to that of $r_1$, with the same parameter $\alpha$ fixed.\\
For all the three cases listed above, before the segmentation step we filter the image by the linear or nonlinear filter described in Sect. \ref{subsec:filter}. 
We will compare the performance of the four numerical schemes presented in Sect. \ref{subsec:edge-detector}. 
For comparison, we will also show the segmented image obtained by the software \emph{SExtractor} \cite{SExtractor}, one of the most popular software in the astronomical community. In these images, each source is represented by the grey-level obtained as average of the pixels values that compose it. 
Since the images are too big, we will work on smaller images of size $300\times 300$ pixels. Hence, we will have $a=b=30$ and $N=M=300$. 
For the active contour method, we use a rectangular initial front $\gamma_0$, i.e. the external boundary of the image, described as the $0$-level set of  
\begin{equation}
u_0(x,y):=1-\left\lvert\frac{x+y-30}{29.6}\right\rvert -\left\lvert\frac{x-y}{29.6}\right\rvert\,,\quad\forall(x,y)\in\Omega\,,
\end{equation}
visible in Fig. \ref{f160u0}, 
\begin{figure}[h!]
\centering
	\includegraphics[width=0.48\textwidth]{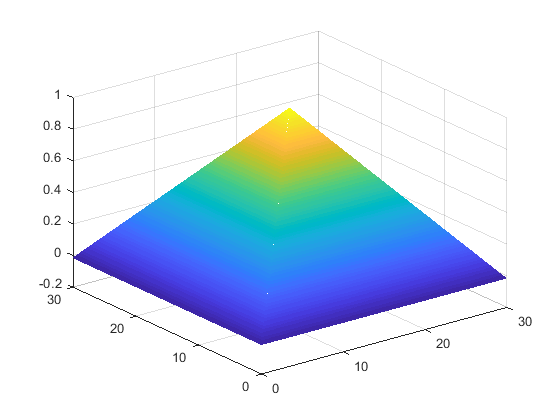}
	\includegraphics[width=0.48\textwidth]{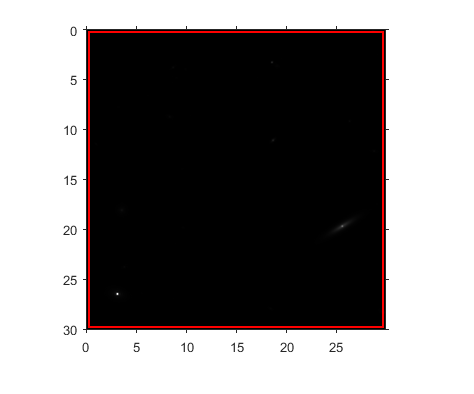}
\caption{Representation function  $u_0$  and related front $\gamma_0$ at time $t=0$.}
\label{f160u0}
\end{figure}
and we filter the image by applying  $5$ iterations of the scheme \eqref{schemacal} 
or by $15$ iterations of the PM method \eqref{eq:PM}, with time step $\widetilde{\Delta t}=10^{-4}$ unless otherwise stated. 
We need to fix two tolerances:  $\varepsilon_F=\Delta x=0.1$, for the identification of the nodes that approximate the front, and $\varepsilon=10^{-3}$ for the stopping criterion. 
The maximum time when the scheme will not converge is set to $T_{max}=50$. 

For the edge-detector function $g_1$, the parameter $p\in\mathbb{N}\setminus\{0\}$ will be fixed according to the contrast in the image between objects and background. If objects  are well defined, we set a value of $p$ small,  if  the edges of the object have pixels with gray tones close to those of the background, the value of $p$ has  to be increased. The function $g_2$ defined in \eqref{g2} in practice does not produce optimal results since it assumes null value only at points where the gradient of the image is maximum and this does not necessarily occur at all points belonging to the edges of the object. Therefore, as proposed in \cite{CFF18}, we modify the function $g_2$, subtracting a constant $c_2 \in [0,1]$ and then rescaling the values in $[0,1]$. The function we use is the following
\begin{equation}\label{g2mod}
\widetilde{g_2}(z):=\frac{1}{1-c_2}\max\{g_2(z)-c_2,0\}\,.
\end{equation}
Thanks to that definition, $\widetilde{g_2}(\lvert\nabla \widetilde{I}_{filt}\rvert$) attains its maximum value equal to $1$ for $\lvert\nabla \widetilde{I}_{filt}\rvert=m$, and null value when $\lvert\nabla \widetilde{I}_{filt}\rvert$ is greater than a fixed threshold, precisely $\lvert\nabla \widetilde{I}_{filt}\rvert\geq(1-c_2)(M-m)+m$. 
In each test,  we provide the values of the parameters involved, e.g. $p$  for the function $g_1$, $c_2$ for the function $\widetilde{g_2}$, and $\nu$ for the dependence from the curvature in the second order schemes  \eqref{DFgk} and \eqref{SLgk}. \\
We acknowledge the contribution of L. Pecci  to the implementation of the methods and to some of the tests presented here. Other numerical experiments are contained in \cite{P18}.

%%% TEST 1
%%%%%%%%%%%%%%%%%%%%%%%%%%%%%%%%%%%%%%%%%%%%%%%%%%%%%%%%%%%%%%%%%%%%%%%%%%%%%
\subsection*{Test 1: \emph{f160.fits}}
The first image, Fig. \ref{f160I0} on the left, is a cropping of a simulated, high resolution astronomical image provided by INAF (Istituto Nazionale di Astrofisica) 
and generated by reproducing data observed by the Hubble Space Telescope (HST). It depicts many stars, galaxies and other celestial bodies, as can be seen from the segmentation obtained with the software \emph{SExtractor}  in Fig. \ref{f160I0} on the right, although it is almost completely black in its original form. 
Our purpose is to apply a segmentation algorithm that traces as many sources as possible, possibly improving the results obtained by \emph{SExtractor} thanks to the introduction of the proposed rescaling functions.

%%%%%%
\subsubsection*{Test1: without rescaling}
Let us start showing the results obtained by the four schemes considered,  without using any rescale function.  
The original image \emph{f160.fits} and the segmentation provided by the software \emph{SExtractor} are shown in Fig. \ref{f160I0}. 
Before applying the active contour schemes, we filter the original image  \emph{f160.fits} by using $5$ iterations of the scheme \eqref{schemacal} with a time step $\widetilde{\Delta t}=10^{-4}$. 
All the active contour methods only identify the brightest celestial body or a few other objects.  
The results are shown in Figs. \ref{f160DFg2} - \ref{f160SLg1k}, the values of the parameters used are mentioned in the captions. 
We only show the results obtained by the schemes \eqref{DFg} and \eqref{SLg} with edge-stopping function $\widetilde{g_2}$ (Figs. \ref{f160DFg2}-\ref{f160SLg2}) and the second order schemes \eqref{DFgk} and \eqref{SLgk} with function $g_1$ (Figs. \ref{f160DFg1k}-\ref{f160SLg1k}). 
Even if we increase the values of the parameters $p$ and $c_2$, we do not get better results. Due to the similarity, we decide to omit the results obtained by using the first order schemes \eqref{DFg} and \eqref{SLg} with edge-detector function $g_1$.  \\
Note that by applying the PM  method \eqref{eq:PM} to the original image, we do not get an improvement as shown in  Fig. \ref{f160SLg2_PM}. 
It is important to note that  the results without a rescaling preprocessing are really bad for all the schemes, even if we apply a nonlinear filtering algorithm. 

\begin{figure}[h!]
\centering
\includegraphics[width=0.4\textwidth]{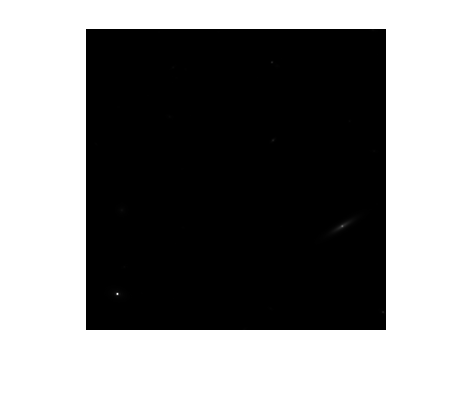}
\includegraphics[width=0.4\textwidth]{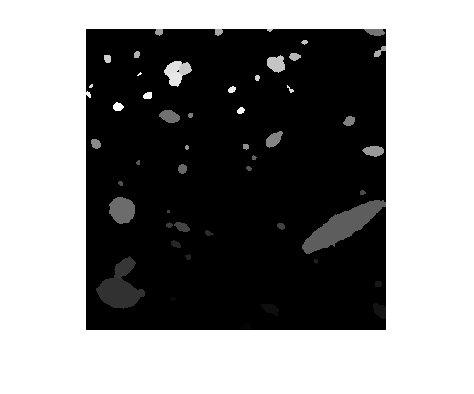}
\vspace{-0.5cm}\caption{Test 1. From left to right: Image \emph{f160.fits}, segmentation of the image provided by the software \emph{SExtractor}.}
\label{f160I0}
\end{figure}

\begin{figure}[h!]
\centering
	\includegraphics[width=0.4\textwidth]{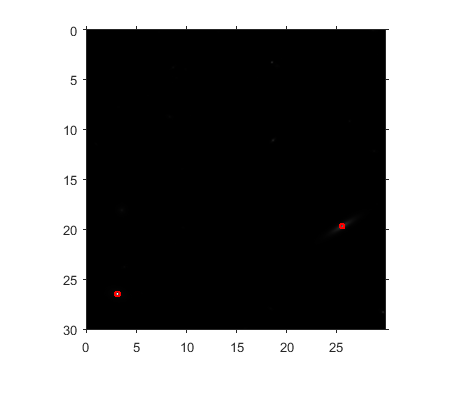}
	\includegraphics[width=0.4\textwidth]{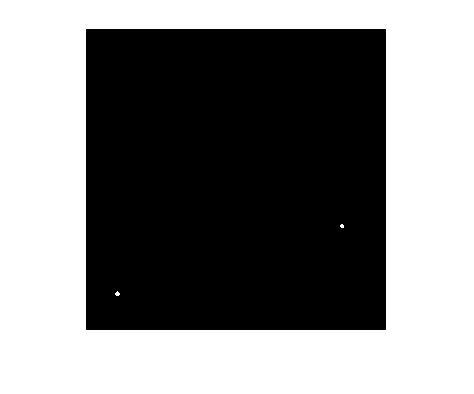}
\vspace{-0.5cm}\caption{Test 1 without rescaling: Position of the front at time $T=15.85$ and segmented image, for the FD scheme \eqref{DFg} by using the edge-detector function $\widetilde{g_2}$, with $c_2=0.8$.} 
\label{f160DFg2}
\end{figure}

\begin{figure}[h!]
\centering
	\includegraphics[width=0.4\textwidth]{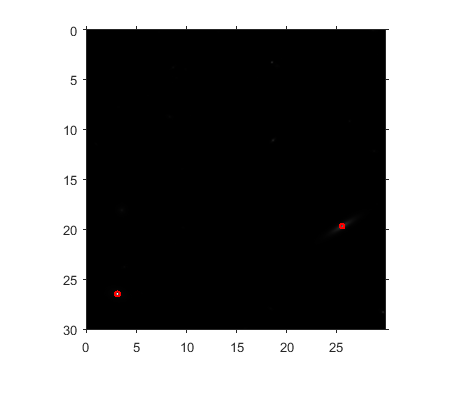}
	\includegraphics[width=0.4\textwidth]{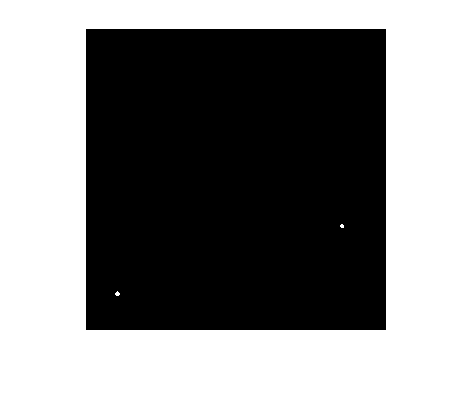}
\vspace{-0.5cm}\caption{Test 1 without rescaling: Position of the front at time $T=15$ and segmented image, for the SL scheme \eqref{SLg} by using the edge-detector function $\widetilde{g_2}$, with $c_2=0.8$.} 
\label{f160SLg2}
\end{figure}

\begin{figure}[h!]
\centering
	\includegraphics[width=0.4\textwidth]{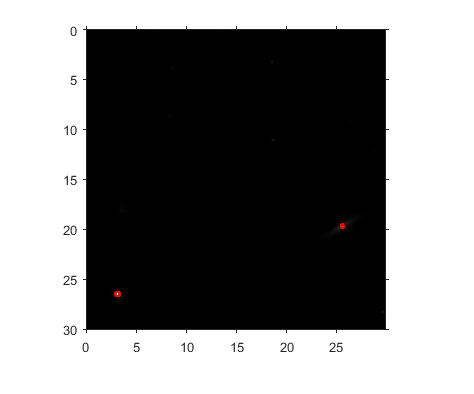}
	\includegraphics[width=0.4\textwidth]{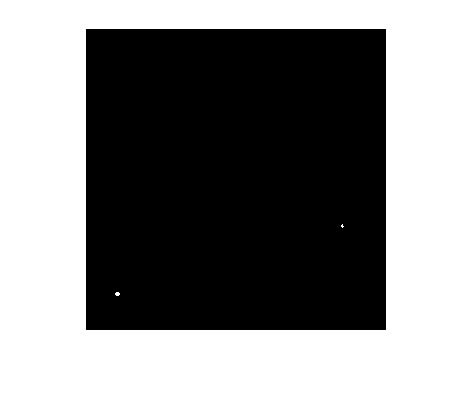}
\vspace{-0.5cm}\caption{Test 1 without rescaling: Position of the front at time  $T=16.61$ and segmented image, for the FD scheme \eqref{DFgk}, by using the edge-detector function  $g_1$, $p=5000$ and $\nu=10^{-4}$.} 
\label{f160DFg1k}
\end{figure}

\begin{figure}[h!]
\centering
	\includegraphics[width=0.4\textwidth]{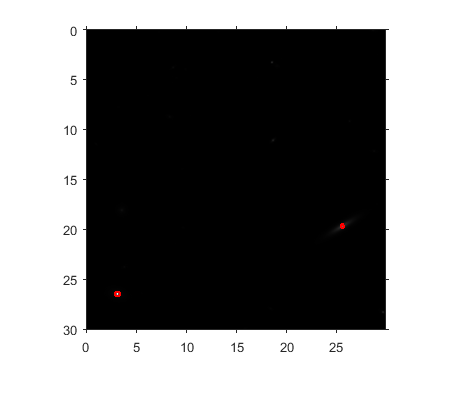}
	\includegraphics[width=0.4\textwidth]{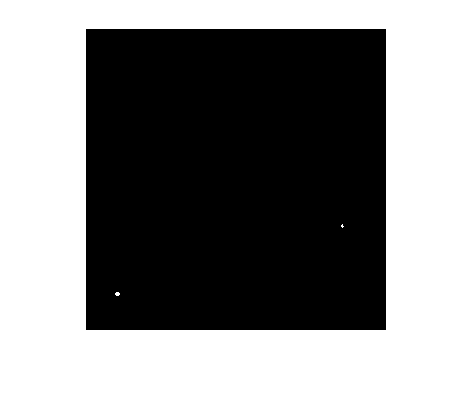}
\vspace{-0.5cm}\caption{Test 1 without rescaling: Position of the front at time  $T=15.4$ and segmented image, for the SL scheme \eqref{SLgk}, by using the edge-detector function $g_1$, $p=5000$ and $\nu=10^{-4}$.}
\label{f160SLg1k}
\end{figure}

\begin{figure}[h!]
	\centering
	\includegraphics[width=0.42\textwidth]{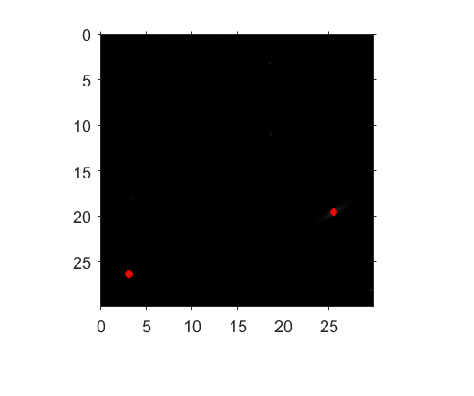}
	\includegraphics[width=0.42\textwidth]{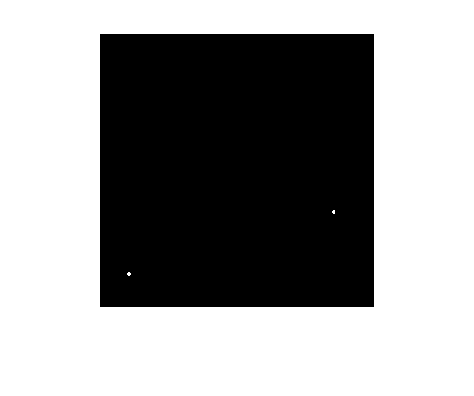}
	\vspace{-0.5cm}\caption{Test 1 without rescaling: Position of the front at time  $T=15$ and segmented image, for the SL scheme \eqref{SLg} by using the edge-detector function $\widetilde{g_2}$, with $c_2=0.8$. Image filtered   by $15$ iterations of the PM method with $f_2$, for $\mu=30$, and  $\widetilde{\Delta t}= 10^{-4}$.}
	\label{f160SLg2_PM}
\end{figure}

%%%%%%%%%%%%%%%%%%%%%%%%%%%%%%%%%%
\subsubsection*{Test 1: rescaling by  $r_1$}
Let us  test the algorithm rescaling the gray tones of the image before applying the segmentation methods. We use $r_1$ setting  $\alpha = 0.25$ since for values of $\alpha$ too close to $1$, the objects are not quite evident, whereas for smaller values the background tones  are amplified excessively. The image  obtained by the rescaling, shown in Fig. \ref{f160_f1I0}, is segmented via the four schemes listed before. 
Also in this case, we omit to show the results obtained by using the first order schemes with edge-detector function $g_1$ since this function fails even if we use a second order scheme, as we can see looking at Figs. \ref{f160_f1DFg1k}-\ref{f160_f1SLg1k}. In fact, $g_1$ does not identify the boundaries, even for large values of the parameter $p$, so that only very few objects are detected. 
Instead, the edge-detector function $\widetilde{g_2}$ (Figs. \ref{f160_f1DFg2}-\ref{f160_f1SLg2}) is able to identify a greater number of objects, even if the approximation of their contours is still non very accurate (for example for the larger galaxy, placed on the right of the image). 
Using the PM nonlinear filtering method after the rescaling process, we can note (comparing Fig. \ref{f160_f1SLg2_PM} and  Fig. \ref{f160_f1SLg2}) that a better segmentation is obtained. In fact, more small objects are detected and visible in the final front and the associated segmented image: see e.g.  two small red points in the central-bottom part of the final front in Fig. \ref{f160_f1SLg2_PM} not present in Fig. \ref{f160_f1SLg2}. 

\begin{figure}[h!]
\centering
\includegraphics[width=0.4\textwidth]{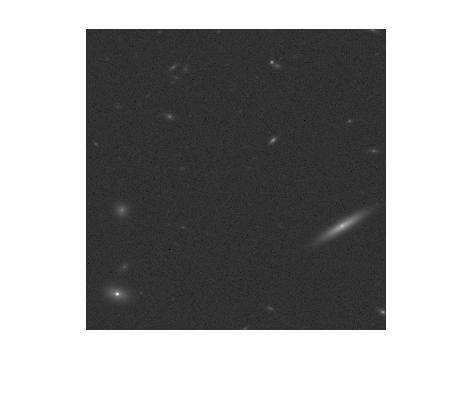}
\includegraphics[width=0.4\textwidth]{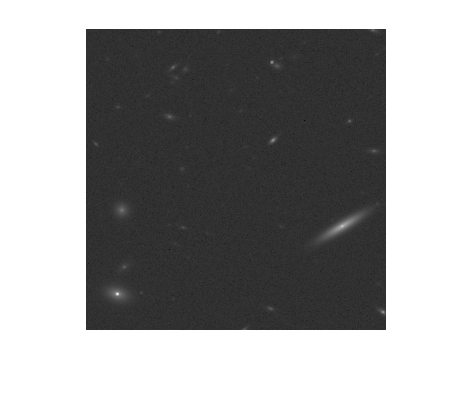}
\vspace{-0.5cm}\caption{Test 1. From left to right: Rescaling of the image \emph{f160.fits} by using the function $r_1$, with $\alpha=0.25$ and its filtered version obtained by  $5$ iterations of the scheme \eqref{schemacal} with time step $\widetilde{\Delta t}=10^{-4}$.}
\label{f160_f1I0}
\end{figure}

\begin{figure}[h!]
\centering
	\includegraphics[width=0.4\textwidth]{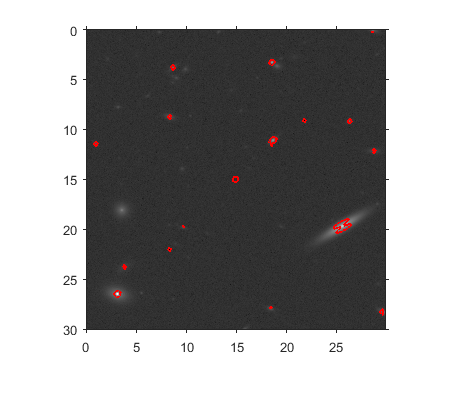}
	\includegraphics[width=0.4\textwidth]{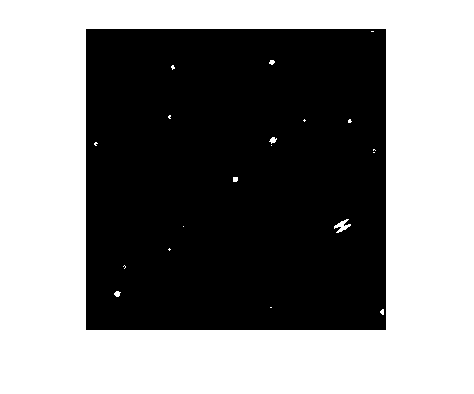}
\vspace{-0.5cm}\caption{Test 1, rescaling by $r_1$: Position of the front at time $T=18.075$ and segmented image, for the FD scheme \eqref{DFg} by using the edge-detector function $\widetilde{g_2}$, with constant $c_2=0.8$.} 
\label{f160_f1DFg2}
\end{figure}

\begin{figure}[h!]
\centering
	\includegraphics[width=0.4\textwidth]{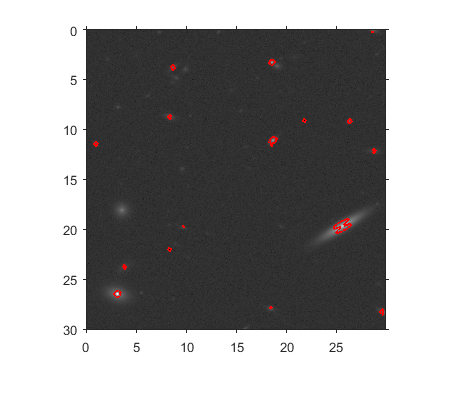}
	\includegraphics[width=0.4\textwidth]{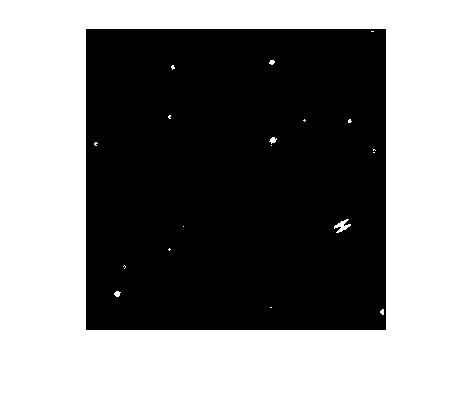}
\vspace{-0.5cm}\caption{Test 1, rescaling by $r_1$: Position of the front at time $T=18.1$ and segmented image, for the SL scheme \eqref{SLg} by using the edge-detector function $\widetilde{g_2}$, with $c_2=0.8$ filtered by the Gaussian filter.} 
\label{f160_f1SLg2}
\end{figure}

\begin{figure}[h!]
	\centering
	\includegraphics[width=0.4\textwidth]{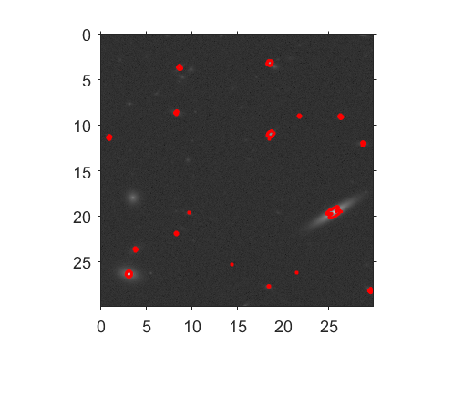}
	\includegraphics[width=0.4\textwidth]{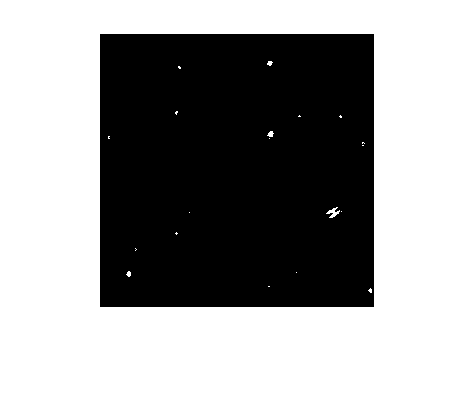}
	\vspace{-0.5cm}\caption{Test 1, rescaling by $r_1$: Position of the front at time $T=18.6$ and segmented image, for the SL scheme \eqref{SLg} by using the edge-detector function $\widetilde{g_2}$, with $c_2=0.8$. The rescaled image has been filtered by $15$ iterations of the PM method with $f_2$, for $\mu=30$, and $\widetilde{\Delta t}= 10^{-4}$.}
	\label{f160_f1SLg2_PM}
\end{figure}

\begin{figure}[h!]
\centering
	\includegraphics[width=0.4\textwidth]{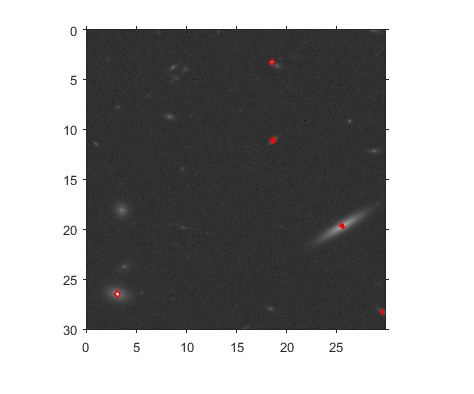}
	\includegraphics[width=0.4\textwidth]{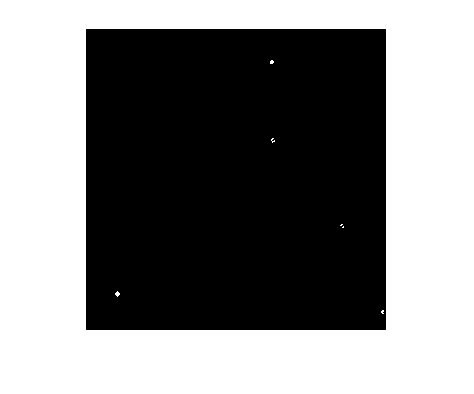}
\vspace{-0.5cm}\caption{Test 1, rescaling by $r_1$: Position of the front at time $T=16.08$ and segmented image, for the FD scheme \eqref{DFgk} by using the edge-detector function $g_1$, $p=5000$ and $\nu=10^{-6}$.} 
\label{f160_f1DFg1k}
\end{figure}

\begin{figure}[h!]
\centering
	\includegraphics[width=0.4\textwidth]{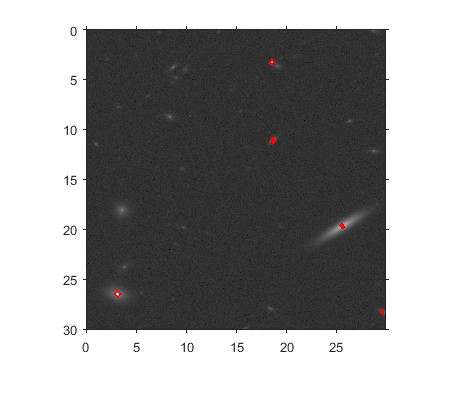}
	\includegraphics[width=0.4\textwidth]{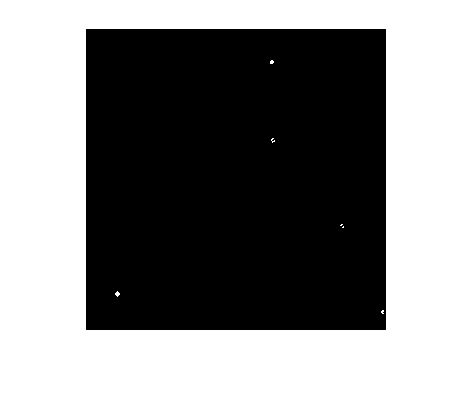}
\vspace{-0.5cm}\caption{Test 1, rescaling by $r_1$: Position of the front at time $T=15.1$ and segmented image, for the SL scheme \eqref{SLgk}  by using the edge-detector function $g_1$ and $\nu=10^{-6}$.} 
\label{f160_f1SLg1k}
\end{figure}

%%%%%%%%%%%%%%%%%%%%%%%%%%%%%%%%%%%%%%%%%%%%%%%%%%%%%%%%%%%%%%%%%%%%%%%%%
\subsubsection*{Test 1: rescaling by $r_3$}
Finally, we present the results obtained by $r_3$, this case seems to give the  best results. 
The parameter chosen is $\beta=8$, for which the boundaries of the objects appear more evident, with tones distant from those of the background (see Fig. \ref{f160_f3I0}). In this case, all the schemes seem to provide satisfactory results, even those based on the use of the edge-detector function $g_1$. 
Due to their similarity, also in this can we show only the two second order schemes with edge-detector function $g_1$. The resulting segmentations are illustrated in Figs. \ref{f160_f3DFg2}-\ref{f160_f3SLg1k}. 
These results show very well the improvements obtained by the rescaling $r_3$. 

\begin{figure}[h!]
\centering
\includegraphics[width=0.4\textwidth]{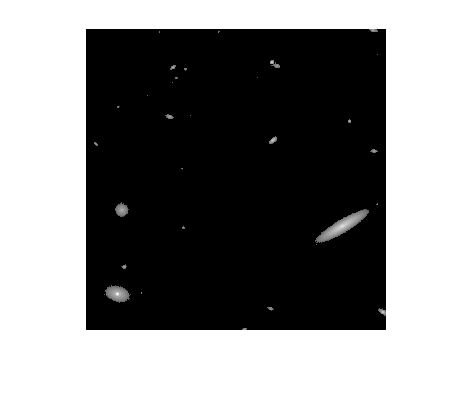}
\includegraphics[width=0.4\textwidth]{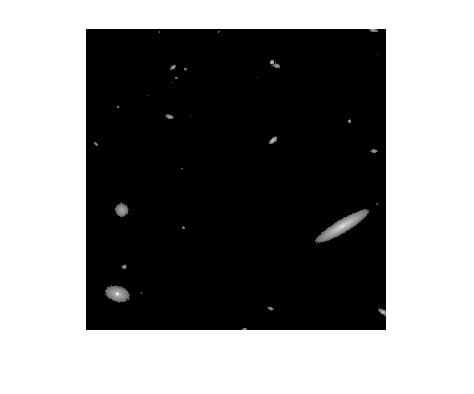}
\vspace{-0.5cm}\caption{Test 1. From left to right: Rescaling of the image \emph{f160.fits} by using the function $r_3$, with $\beta=8$, and its filtered version obtained by 
 $5$ iterations of the scheme \eqref{schemacal} with time step $\widetilde{\Delta t}=10^{-4}$.}
\label{f160_f3I0}
\end{figure}

\begin{figure}[h!]
\centering
	\includegraphics[width=0.4\textwidth]{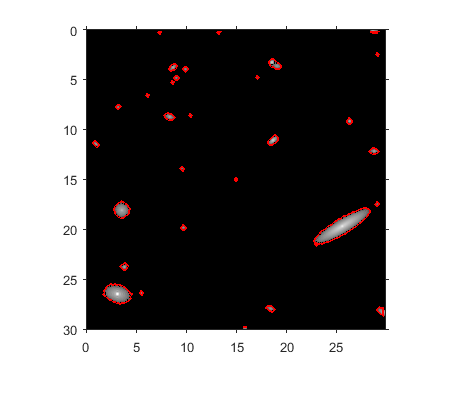}
	\includegraphics[width=0.4\textwidth]{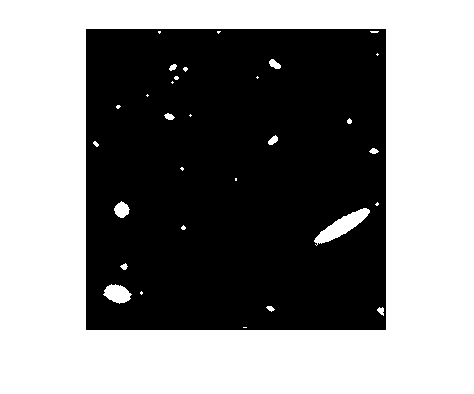}
\vspace{-0.5cm}\caption{Test 1, rescaling by $r_3$: Position of the front at time $T=15.4$ and segmented image, for the FD scheme \eqref{DFg} by using the edge-detector function $\widetilde{g_2}$, with $c_2=0.8$.} 
\label{f160_f3DFg2}
\end{figure}

\begin{figure}[h!]
\centering
	\includegraphics[width=0.4\textwidth]{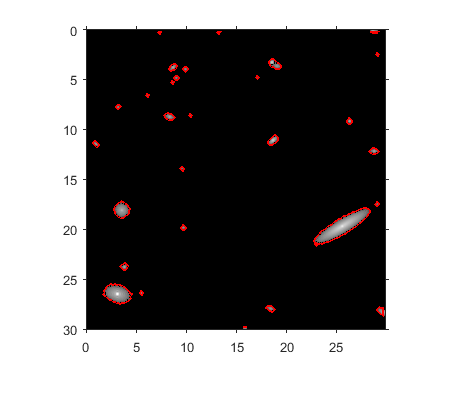}
	\includegraphics[width=0.4\textwidth]{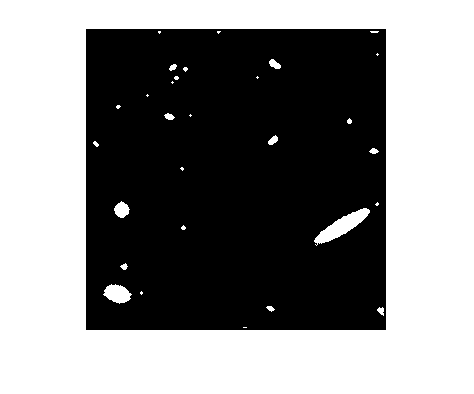}
\vspace{-0.5cm}\caption{Test 1, rescaling by $r_3$: Position of the front at time $T=14.9$ and segmented image, for the SL scheme \eqref{SLg} by using the edge-detector function $\widetilde{g_2}$, with constant $c_2=0.8$.} 
\label{f160_f3SLg2}
\end{figure}

\begin{figure}[h!]
\centering
	\includegraphics[width=0.4\textwidth]{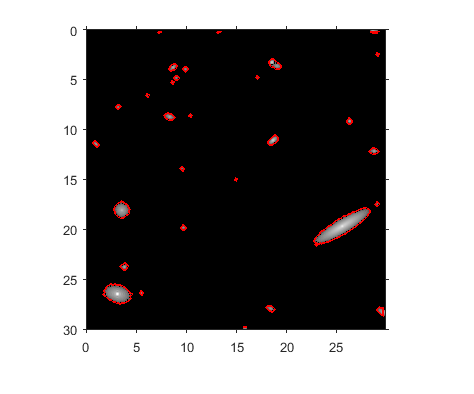}
	\includegraphics[width=0.4\textwidth]{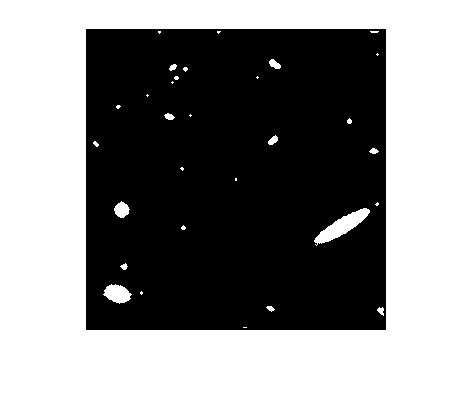}
\vspace{-0.5cm}\caption{Test 1, rescaling by $r_3$: Position of the front at time $T=15.55$ and segmented image, for the FD scheme \eqref{DFgk} by using the edge-detector function $g_1$, $p=5000$ and $\nu=10^{-4}$.} 
\label{f160_f3DFg1k}
\end{figure}

\begin{figure}[h!]
\centering
	\includegraphics[width=0.4\textwidth]{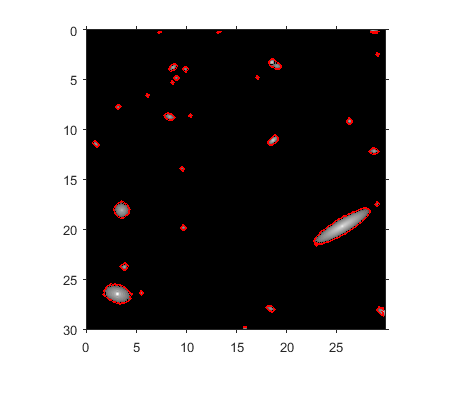}
	\includegraphics[width=0.4\textwidth]{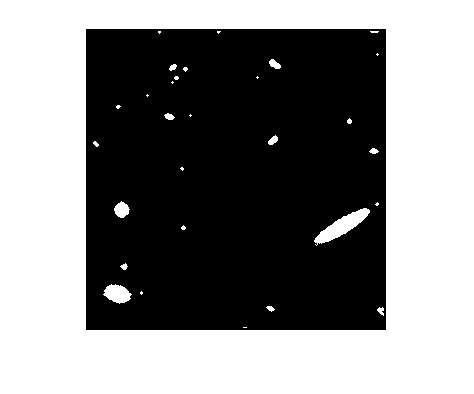}
\vspace{-0.5cm}\caption{Test 1, rescaling by $r_3$: Position of the front at time $T=28.8$ and segmented image, for the SL scheme \eqref{SLgk} by using the edge-detector function $g_1$, $p=5000$ and $\nu=10^{-4}$.} 
\label{f160_f3SLg1k}
\end{figure}

%%%%%%%%%%%%%%%%%%%%%%%%%%%%%%%%%%%%%%%%%%%%%%%%%%%%%%%%%%%%%%%%%%%%%%%%%%%%%
\subsection*{Test 2: \emph{real.fits}}
We now consider a clipping of a real low resolution image  generated by the Hubble Space Telescope and provided by INAF.  
This image has been acquired by observing a portion of the sky at high depth, in order to identify a very large number of sources. 
However, this technique leads to an increase in the amount of noise present in the image, as can be seen looking at the left image in Fig. \ref{realI0&regI0}. \\
Let us compare the performances of  different schemes with or without a rescaling process. The input images we consider for the segmentation algorithms are reported in Fig. \ref{realI0&regI0}. On the left we can see the original image that we store in a file called \emph{real.fits}, in the middle we find the image obtained by using the rescaling function $r_1$, and the analogous obtained rescaling by $r_3$ (on the right). 
Due to the high level of noise, here we increase the time step $\widetilde{\Delta t}$ (from $10^{-4}$ to $10^{-3}$) in the filtering process to obtain the corresponding filtered images. Then we use the rescaling on the filtered images. 
As can be easily noted looking at Fig. \ref{realI0&regI0}, both the proposed rescaling functions improve a lot  the visibility of the celestial objects present in the  image. The resulting image obtained by  $r_3$ seem to be better. 
\begin{figure}[t]
\hspace{-0.9cm}
\centering
\includegraphics[width=0.37\textwidth]{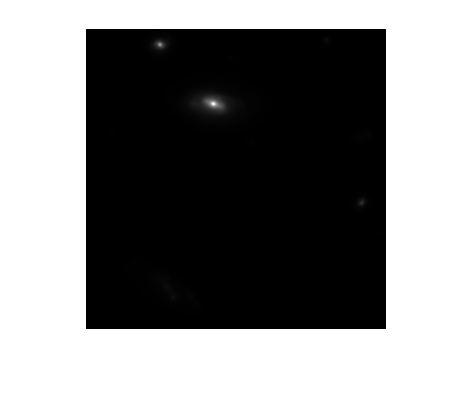}     	
\hspace{-0.9cm}
\includegraphics[width=0.37\textwidth]{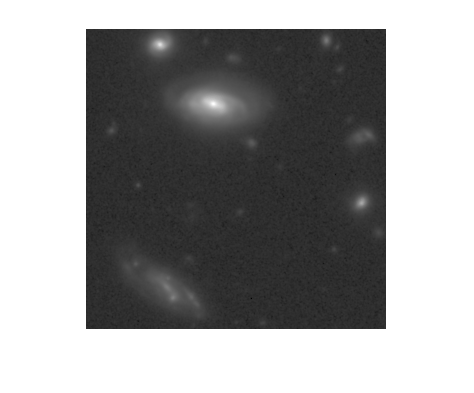}		 \hspace{-0.9cm}
\includegraphics[width=0.37\textwidth]{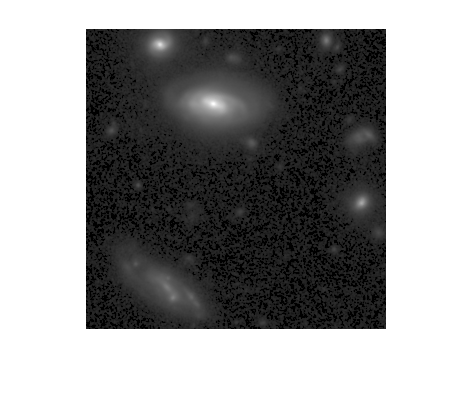}		
\caption{Test 2. From left to right: Original image \emph{real.fits}, rescaling of the image \emph{real.fits} by using the function $r_1$ with $\alpha=0.25$, rescaling of the image \emph{real.fits} by using the function $r_3$, with $\beta=4$.}
\label{realI0&regI0}
\end{figure}
\begin{figure}[t]
\centering
\includegraphics[width=0.4\textwidth]{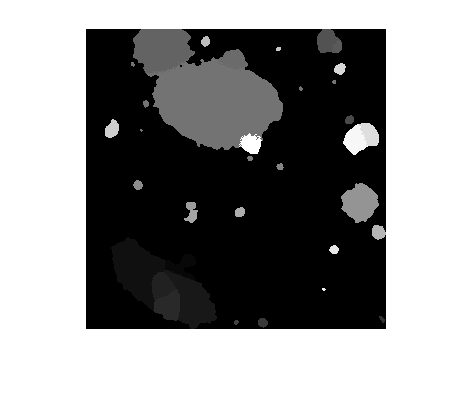}
\vspace{-0.5cm}
\caption{Test 2. Segmentation of the image \emph{real.fits} provided by the software \emph{SExtractor}.}
\label{realI0SExt}
\end{figure}

Let us start commenting the segmentation results obtained by the different schemes without any rescaling process. 
In Fig. \ref{realI0SExt}, we report the segmentation of the image obtained by applying the software \emph{SExtractor}. 
In Fig. \ref{real:segmNoRisc} we can see the performances of the SL scheme \eqref{SLg} for two different choices of the edge-detector function ($g_1$ and $\widetilde{g_2}$) and the SL scheme \eqref{SLgk} with edge-detector function $g_1$. We report only the results for SL schemes since by FD schemes we obtained very similar results. Note  that all the different schemes recognize only the two objects visible in the original image reported on the top-left of Fig. \ref{realI0&regI0}, 
so they are far away from the real configuration of celestial bodies. 

\begin{figure}[h!]
\centering
	\includegraphics[width=0.4\textwidth]{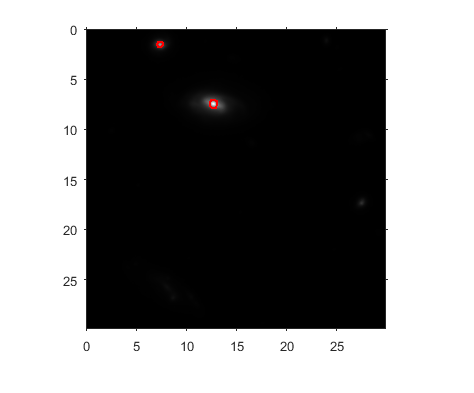}		
	\includegraphics[width=0.4\textwidth]{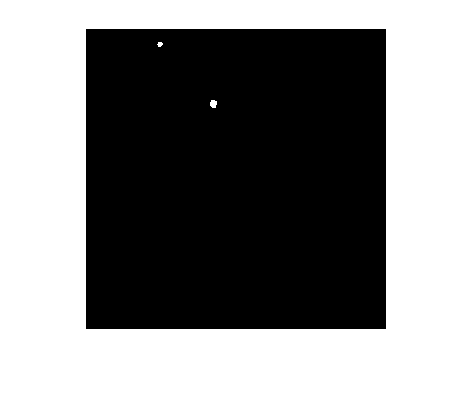}\\
\vspace{-1cm}
	\includegraphics[width=0.4\textwidth]{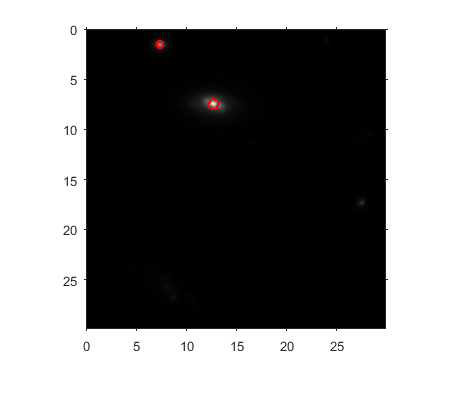} 		
	\includegraphics[width=0.4\textwidth]{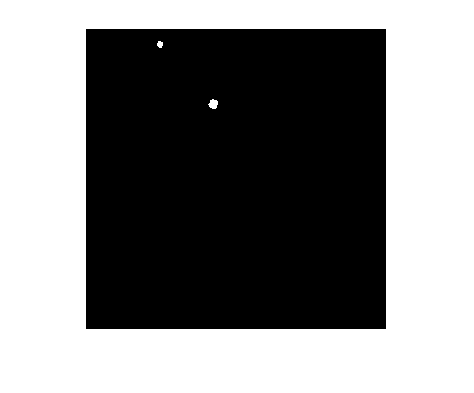}\\	
\vspace{-1cm}
	\includegraphics[width=0.4\textwidth]{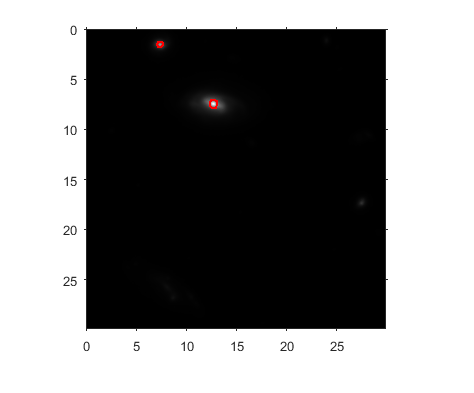}	
	\includegraphics[width=0.4\textwidth]{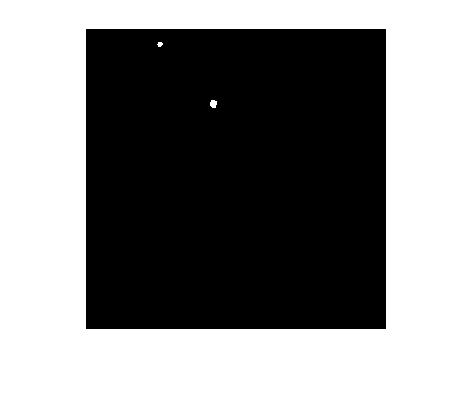}
\vspace{-0.5cm}
\caption{Test 2 without rescaling. From top to bottom: Position of the front and segmented image for the SL scheme \eqref{SLg} by using the edge-detector function $g_1$, with  $p=10^4$. 
Same scheme by using edge-detector function $\widetilde{g_2}$, with  $c_2=0.6$. 
 Position of the front and segmented image for the SL scheme  \eqref{SLgk} by using the edge-detector function $g_1$, with  $p=10^4$ and $\nu=10^{-4}$.}
\label{real:segmNoRisc}
\end{figure}
Therefore we need a rescaling process to improve the results. 
Looking at the results obtained by $r_1$, we note that the number of detected objects is highly  improved. We report in Fig. \ref{real:segm-f1} the results obtained by the schemes FD and SL  only with edge-detector $\widetilde{g_2}$, using the edge-function $g_1$ the front collapses until it disappears from the figure. This is due to the fact that the edges of the objects are very blurred and, even if we choose high values for the parameter $p$, the variations of gray tones do not allow to detect the presence of an object. 
On the contrary, using $\widetilde{g_2}$ we can find  an adequate number of objects, but we cannot detect accurately  the boundaries of many galaxies (see Fig. \ref{real:segm-f1}). Anyway, this result is more accurate than the performance provided by the software \emph{SExtractor} (Cf. Fig. \ref{realI0SExt}). 
Looking more in details Fig. \ref{real:segm-f1}, some differences between the FD and SL schemes are visible, even if are small (e.g. at the bottom-right part of the big central-upper galaxy we can note a connected part for the FD scheme, which is splitted  by SL scheme). 
For comparison reasons, we report in Fig. \ref{real:segm-f1_PM} the result obtained by the same FD scheme with edge-detector function $\widetilde{g_2}$, but the image rescaled by the function $r_1$ is filtered by the PM method before applying the FD segmentation scheme. The position of the final front and the segmented image reported in Fig. \ref{real:segm-f1_PM} show that, applying a nonlinear filtering algorithm as the PM method before the segmentation process, the results can be improved (a lot of small stars are recognized), but a rescaling process is still necessary even if we apply that filtering method. 
\begin{figure}[h!]
\centering
	\includegraphics[width=0.4\textwidth]{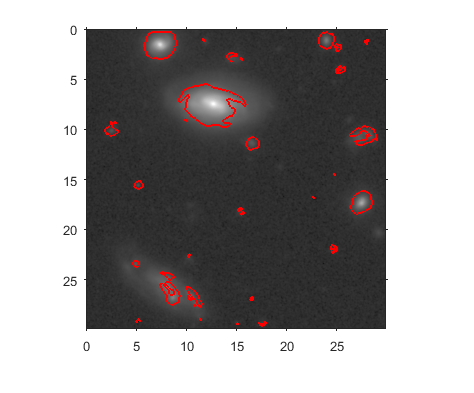}	
	\includegraphics[width=0.4\textwidth]{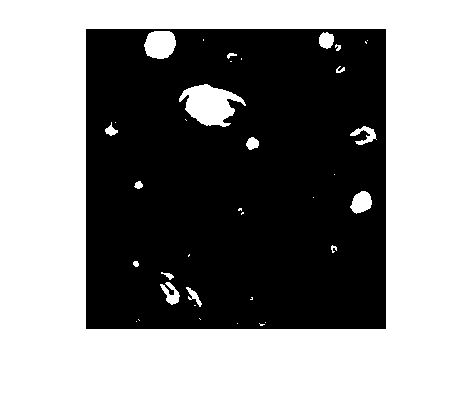} \\
\vspace{-0.4cm}
	\includegraphics[width=0.4\textwidth]{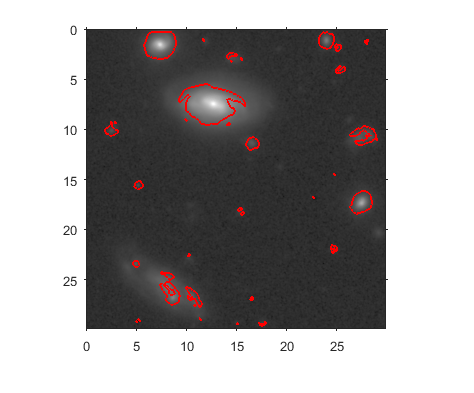}  
	\includegraphics[width=0.4\textwidth]{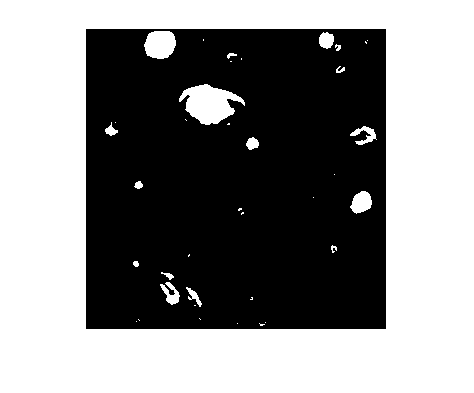}
\vspace{-0.5cm}
\caption{Test 2, rescaling by $r_1$. Position of the front and segmented image for the FD scheme \eqref{DFg} (first row) and the SL scheme \eqref{SLg} (second row) by using the edge-detector function $\widetilde{g_2}$, with  $c_2=0.78$.} 
\label{real:segm-f1}
\end{figure}

\begin{figure}[h!]
	\centering
	\includegraphics[width=0.44\textwidth]{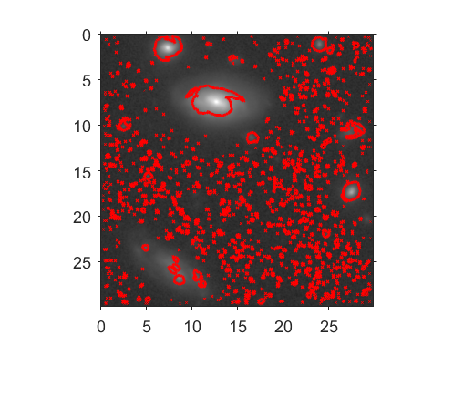}		
	\includegraphics[width=0.44\textwidth]{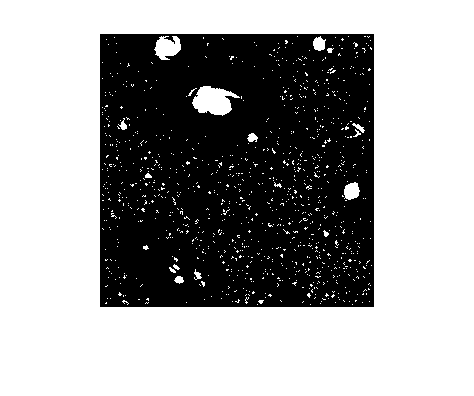} \\
	\vspace{-0.5cm}
	\caption{Test 2, rescaling by $r_1$, filtered by $15$ iterations of the PM method with $f_2$, for $\mu=30$, and $\widetilde{\Delta t} = 10^{-3}$. Position of the front and segmented image for the FD scheme \eqref{DFg} by using the edge-detector function $\widetilde{g_2}$, with  $c_2=0.78$.} 
	\label{real:segm-f1_PM}
\end{figure}

Finally, let us analyze the results obtained by applying the rescaling function $r_3$. The parameters chosen in that case is $\beta=4$ 
(see the right image in Fig. \ref{realI0&regI0}), since greater values provide apparently worst quality. This is due to the poor performance of the Otsu method in this case,  note that this method identifies false sources among the pixels of the background. For the type of results provided by the different active contours, similar observations apply to the images obtained by the  rescaling $r_1$, as we can observe from the segmentations shown in Fig. \ref{real:segm-f3}. With the $g_1$ function, the front collapses on itself, disappearing without identifying any object. 
The results obtained by using $\widetilde{g_2}$ are better even if not satisfactory, due to the noise component, which is excessively amplified, as visible in Fig. \ref{real:segm-f3}. 
\begin{figure}[h!]
\centering
	\includegraphics[width=0.4\textwidth]{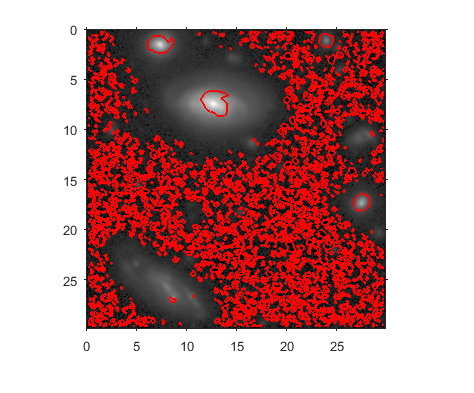}		
	\includegraphics[width=0.4\textwidth]{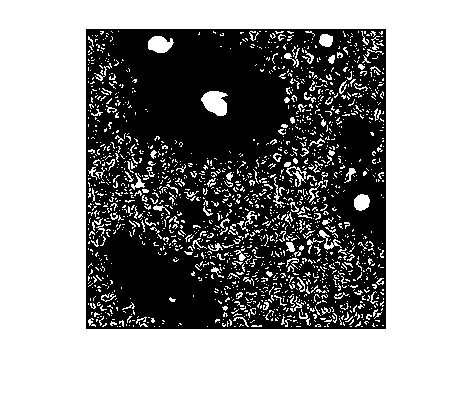} \\
\vspace{-0.4cm}
	\includegraphics[width=0.4\textwidth]{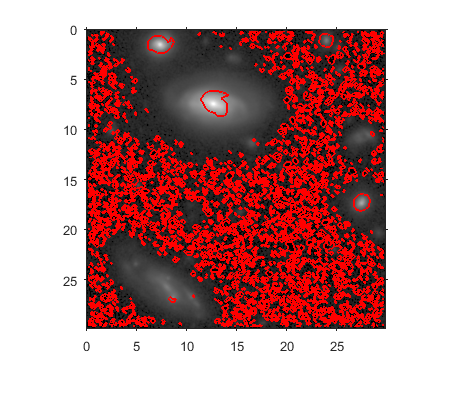}	
	\includegraphics[width=0.4\textwidth]{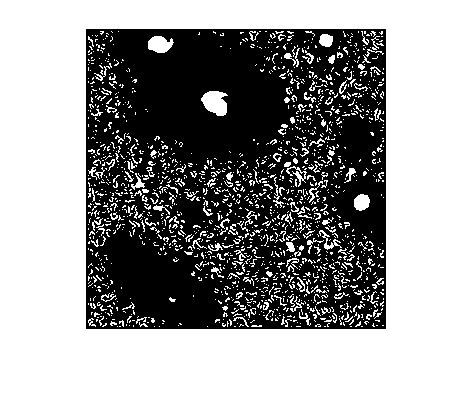}
\vspace{-0.5cm}
\caption{Test 2, rescaling by $r_3$. Position of the front and segmented image for the FD scheme \eqref{DFg} (first row) and the SL scheme \eqref{SLg} (second row) by using the edge-detector function $\widetilde{g_2}$, with  $c_2=0.6$.} 
\label{real:segm-f3}
\end{figure}

%%%%%%%%%%%%%%%%%%%%%%%%%%%%%%%%%%%%%%%%%%%%%%%%%%%%%%%%%%%%%%%%%%%%%%%%%%
\section{Conclusions and future perspectives}\label{sec:conclusions}
We have proposed different rescaling functions in order to improve the detection of objects in astronomical images, identifying  a greater number of celestial bodies. In particular, the use of the function $r_3$ has  improved a lot  the visibility, getting better results, in particular for high-resolution images as \emph{f160.fits}. Unfortunately, when the SNR is very low, the results are still not satisfactory, although we notice an improvement with respect to the solutions provided by classical methods without rescaling. This failure is due to the inaccurate threshold selected via the Otsu algorithm before applying the rescaling by using the function $r_3$. For  low-resolution images, the function $r_1$ seems to provide the best results, even compared with those produced by the software \emph{SExtractor} commonly used by astronomers. 
Future improvements of the method can focus on different threshold algorithms to select the optimal one used by the rescaling function $r_3$, hopefully this will allow for a correct classification of the pixels belonging to the background. 
We also considered a filtering pre-processing step before segmenting, comparing the linear Gaussian filter and the nonlinear PM method.
What we noted is that the PM nonlinear method improves the results detecting few more objects, but a rescaling preprocessing is necessary also in this case, the segmentation fails without it. 
We have compared the performances of first and second order FD and SL schemes, using different parameters and two different edge-detector functions. 
From the numerical  simulations on virtual and real images, we can conclude that the edge-stopping function $g_1$ is not a good choice. In fact, the light sources often have not well defined outlines, so that this function can not correctly identify them. The edge-detector function $\widetilde{g_2}$ provides the best results,  and is able to detect a greater number of celestial objects.  However, these methods are still not optimal in the case of very disturbed images. In the future, we want to explore different methods, as for example high-order ``filtered'' schemes recently proposed \cite{FPT18,FPT18a,FPT18b} or the active contour without edges scheme proposed  by Chan and Vese in \cite{CV01, VC02}, able (perhaps) to better identify objects with blurred  and not well defined edges. 
Moreover, we would like to analyze in more detail the performances of other filtering methods, in order to find an appropriate choice to reduce the huge amount of noise that is a typical feature of astronomical images. Some attempts in this direction are shown in \cite{RTCMFF}, they confirm that this will be a difficult task.

\paragraph*{Acknowledgements}
We would like to thank the National Group INdAM-GNCS for the financial support given to this research and the Istituto Nazionale di Astrofisica placed in Rome for the input data. 
This research has been carried on within the INdAM-INAF project FOE 2015 ``OTTICA ADATTIVA''.

\bibliographystyle{plain}
\bibliography{referenc}

\end{document}